\documentclass [12pt] {amsart}

\usepackage{amsthm}
\usepackage{amssymb}

\usepackage{amsmath}
\usepackage{epsfig}

\newtheorem{theorem}{Theorem}

\newtheorem{definition}[theorem]{Definition}

\newtheorem{lemma}[theorem]{Lemma}

\begin{document}

\title{The Rotation Class of a Flow}

\author{Alex Clark}

\address{Department of Mathematics, University of North Texas, 
Denton, TX 76203-1430}

\email{alexc@unt.edu}

\date{October 31, 2000.}

\subjclass{37B05,37B45}

\thanks{This work was funded in part by a faculty
research grant from the University of North Texas.}
\maketitle

\begin{abstract}
Generalizing a \ construction of A. Weil, we introduce a topological
invariant for flows on compact, connected, finite dimensional, abelian,
topological groups. We calculate this invariant for some examples and
compare the invariant with other flow invariants.
\end{abstract}

\section{Introduction}

For our purposes a \emph{flow }is a continuous group action of $\left( 
\mathbf{R},+\right) $. We consider the flow $\phi $ on $X$ and the flow $
\psi $ on $Y$ to be \emph{topologically equivalent }when there is a
homeomorphism $h:X\rightarrow Y$ which takes orbits of $\phi $ onto orbits
of $\psi $ in such a way that the orientation of orbits is preserved, in
which case we write $h:\phi \overset{\text{top}}{\approx }\psi .$ Following
the lead of A. Weil as outlined in \cite{W1}\cite{W2}, we introduce the 
\emph{rotation class }of a flow: an invariant of topological equivalence for
flows on finite dimensional, compact, connected, abelian, topological groups
-- hereafter referred to as $n-$solenoids, see definition \ref{sol} or \cite
{C2}. A modern account of Weil's torus invariant and its application to
toral flows may be found in \cite{AZ},\cite[Chapter 6]{ABZ}. For each $n-$
solenoid we shall describe a covering space and a compactification of the
covering space, the \emph{compactification of perspective}. The remainder of
this compactification is homeomorphic to $S^{n-1}.$ Informally stated, we
shall show that if $\phi $ and $\psi $ are topologically equivalent flows on
the $n-$solenoid $\sum_{\overline{M}}$ and if the lifted flow of $\phi $ has
an orbit which has an $\omega -$limit point in the remainder of the
compactification of perspective, then any corresponding orbit in the lifted
flow of $\psi $ will have a corresponding $\omega -$limit point in the
remainder of the compactification. The group structure of $\sum_{\overline{M}
}$ determines which points in the remainder may correspond.

We compare this invariant with other flow invariants from \cite{C1} and \cite
{MZ} and calculate these invariant in some examples. In the course of doing
this we provide a general technique for calculating the exponent group as
introduced in \cite{C1}. It will then be apparent that the rotation class is
well suited for investigating topological equivalence but that it is not as
well suited as other invariants for determining ergodic properties of the
associated flow.

\section{Background and the Compactification of Perspective}

\begin{definition}
\label{sol}For a fixed $n$ and a sequence $\overline{M}=\left(
M_{1},M_{2},...\right) $ of $n\times n$ matrices $M_{i}$ with integer
entries and non-zero determinants, we define the topological group $
\sum\nolimits_{\overline{M}}$ with identity $e_{\overline{M}}$ to be the
inverse limit of the inverse sequence $\{\mathbf{X}_{j},f_{j}^{i}\}$, where $
\mathbf{X}_{j}=\mathbf{T}^{n}$ for all\ $j\in \mathbf{N}$ and $f_{j}^{j+1}$
is the topological epimorphism represented by the matrix $M_{j}$;\ $
f_{j}^{j+1}\circ p^{n}=p^{n}\circ M_{j}$. 
\begin{equation*}
\sum\nolimits_{\overline{M}}\overset{\text{def}}{=}\left\{ \left\langle
x_{j}\right\rangle _{j=1}^{\infty }\in \prod_{j=1}^{\infty }\mathbf{T}
^{n}|f_{j}^{j+1}\left( x_{j+1}\right) =x_{j}\text{ for all }j\in \mathbf{N}
\right\} \text{, }
\end{equation*}
and we define such an inverse limit $\sum\nolimits_{\overline{M}}$ to be an $
n$\emph{-solenoid}.
\end{definition}

We shall only consider the case that $n\in \left\{ 2,3,...\right\} $. As in 
\cite{C2} we have the fibration with unique path lifting 
\begin{eqnarray*}
\pi _{\overline{M}} &:&\mathbf{R}^{n}\rightarrow \sum\nolimits_{\overline{M}}
\\
\pi _{\overline{M}}\left( s\right) &=&\left( \pi ^{n}\left( s\right) ,\pi
^{n}\left( M_{1}^{-1}\left( s\right) \right) ,...,\pi ^{n}\left(
M_{n}^{-1}\circ \cdots \circ M_{1}^{-1}\left( s\right) \right) ,...\right) ,
\end{eqnarray*}
where $\pi ^{n}:\mathbf{R}^{n}\rightarrow \mathbf{T}^{n}=\left( \mathbf{R/Z}
\right) ^{n}$ is the quotient covering map. If $p_{i}$ denotes the
projection of $\sum_{\overline{M}}$ onto the $i^{th}$ $\mathbf{T}^{n}$
factor, then $p_{1}^{-1}\left( e\right) $ is a Cantor set or a finite
discrete space. We may then form the following covering of $\sum\nolimits_{
\overline{M}}.$

\begin{definition}
\begin{equation*}
E_{\overline{M}}\overset{\text{def}}{=}\mathbf{R}^{n}\times p_{1}^{-1}\left(
e\right)
\end{equation*}
and 
\begin{equation*}
\Pi _{\overline{M}}:E_{\overline{M}}\rightarrow \sum\nolimits_{\overline{M}}
\end{equation*}
is the covering map given by 
\begin{equation*}
\Pi _{\overline{M}}\left( s,x\right) =\pi _{\overline{M}}\left( s\right) +x.
\end{equation*}
\end{definition}

See \cite{M} for a similarly defined covering map.\bigskip

With $h_{n}:\mathbf{R}^{n}\rightarrow \mathbf{D}^{n}=\left\{ x\in \mathbf{R}
^{n}|\left\| x\right\| <1\right\} $ denoting the homeomorphism 
\begin{equation*}
x=\left( x_{1},..,x_{n}\right) \mapsto \frac{1}{\sqrt{1+x_{1}^{2}+\cdots
+x_{n}^{2}}}\cdot x
\end{equation*}
and identifying a point $x$ of $\mathbf{R}^{n}$ with its $h_{n}$ image in $
\mathbf{R}^{n}$, we may consider $\overline{\mathbf{D}^{n}}=\left\{ x\in 
\mathbf{R}^{n}|\left\| x\right\| \leq 1\right\} $ to be a compactification
of $\mathbf{R}^{n}.$ Similarly we may consider $\overline{\mathbf{D}^{n}}
\times p_{1}^{-1}\left( e\right) $ to be a compactification of $E_{\overline{
M}}$ via $h_{n}\times id_{p_{1}^{-1}\left( e\right) }$. With $
S^{n-1}=\left\{ x\in \mathbf{R}^{n}|\left\| x\right\| =1\right\} ,$ we form
the equivalence relation $\approx _{\overline{M}}$ on $\overline{\mathbf{D}
^{n}}\times p_{1}^{-1}\left( e\right) $ defined by 
\begin{equation*}
\left[ \left( s,x\right) \approx _{\overline{M}}\left( s^{\prime },x^{\prime
}\right) \right] \iff \left[ s=s^{\prime }\in S^{n-1}\right] .
\end{equation*}
We then have the quotient mapping 
\begin{equation*}
q_{\overline{M}}:\overline{\mathbf{D}^{n}}\times p_{1}^{-1}\left( e\right)
\rightarrow \left( \overline{\mathbf{D}^{n}}\times p_{1}^{-1}\left( e\right)
\right) /\approx _{\overline{M}},
\end{equation*}
the image of which is also a compactification of $E_{\overline{M}},$ \emph{
the compactification of perspective}$.$ The remainder of this
compactification, 
\begin{equation*}
R_{\overline{M}}\overset{\text{def}}{=}\left[ \left( \overline{\mathbf{D}^{n}
}\times p_{1}^{-1}\left( e\right) \right) /\approx _{\overline{M}}\right] -
\left[ q_{\overline{M}}\circ \left( h_{n}\times id_{p_{1}^{-1}\left(
e\right) }\right) \left( E_{\overline{M}}\right) \right]
\end{equation*}
is canonically homeomorphic to $S^{n-1}$
\begin{equation*}
\left[ \left( t,x\right) \right] _{\approx _{\overline{M}}}\rightsquigarrow
t,
\end{equation*}
where $\left[ \left( t,x\right) \right] _{\approx _{\overline{M}}}$ denotes
the $\approx _{\overline{M}}$-class of a point in $R_{\overline{M}}.$
Hereafter we identify $R_{\overline{M}}$ with $S^{n-1}$ as above.

In \cite{C2} the linear flows on $\sum_{\overline{M}}$
\begin{equation*}
\mathcal{F}_{\overline{M}}=\{\Phi _{\overline{M}}^{\mathbf{\omega }}\mid 
\mathbf{\omega }\in \mathbf{R}^{n}-\{\mathbf{0}\}\}\text{ given by}
\end{equation*}
\begin{equation*}
\Phi _{\overline{M}}^{\mathbf{\omega }}:\mathbf{R}\times \sum\nolimits_{
\overline{M}}\rightarrow \sum\nolimits_{\overline{M}};\text{ }\Phi _{
\overline{M}}^{\mathbf{\omega }}\left( t,x\right) =\pi _{\overline{M}}\left(
t\cdot \mathbf{\omega }\right) +x
\end{equation*}
were introduced. Each linear flow $\Phi _{\overline{M}}^{\mathbf{\omega }}$
lifts to the flow $\widetilde{\Phi _{\overline{M}}^{\mathbf{\omega }}}$ on $
E_{\overline{M}}$ given by
\begin{equation*}
\widetilde{\Phi _{\overline{M}}^{\mathbf{\omega }}}\left( t,\left(
s,x\right) \right) =\left( t\cdot \mathbf{\omega }+s\mathbf{,}x\right) .
\end{equation*}
Since $h_{n}\left( t\cdot \mathbf{\omega }+s\right) \rightarrow \dfrac{
\mathbf{\omega }}{\left\| \mathbf{\omega }\right\| }$ as $t\rightarrow
\infty ,$ the image of any $\widetilde{\Phi _{\overline{M}}^{\mathbf{\omega }
}}-$orbit in the compactification of perspective tends asymptotically toward
the same point $\dfrac{\mathbf{\omega }}{\left\| \mathbf{\omega }\right\| }$
in the remainder $S^{n-1}.$ This leads naturally to the following.

\begin{definition}
\label{persp}The point $x\in S^{n-1}$ is a \emph{point of perspective }of
the flow $\phi $ on $\sum\nolimits_{\overline{M}}$ if and only if $\phi $
has an orbit which lifts to an orbit $\mathcal{O}$ in the compactification
of perspective satisfying $\left\{ \mathcal{O}\left( t_{i}\right) \right\}
_{i\in \mathbf{N}}\rightarrow x$ for some sequence of real numbers $\left\{
t_{i}\right\} _{i\in \mathbf{N}}\rightarrow \infty .$
\end{definition}

Thus, the linear flow $\Phi _{\overline{M}}^{\mathbf{\omega }}$ has
precisely one point of perspective: $\dfrac{\mathbf{\omega }}{\left\| 
\mathbf{\omega }\right\| }.$

\center{\epsfig{file=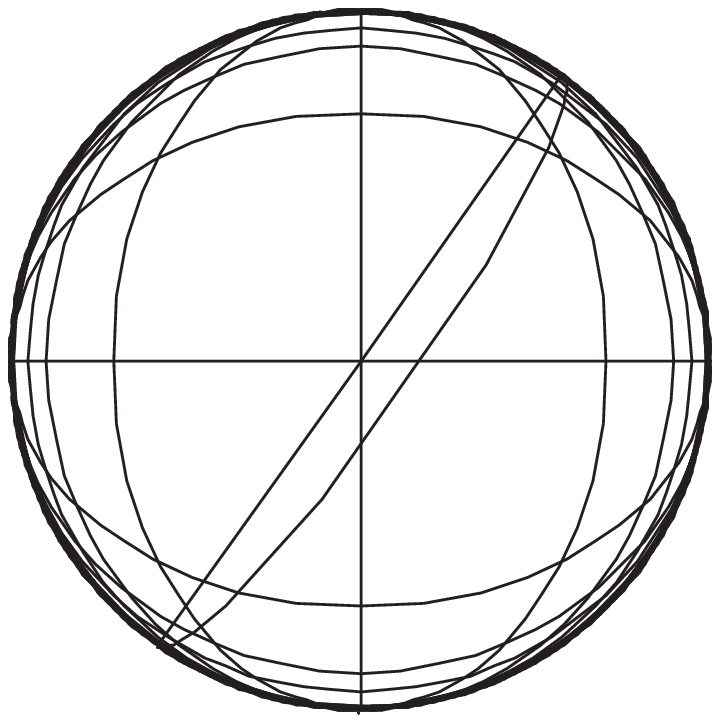}}

\begin{center}
{\small The compactification of perspective of the torus, along with the
integer coordinate grid and two orbits of an irrational flow.}
\end{center}

\section{The Rotation Class of a Flow}

In \cite{C2} the linear flows were classified up to topological equivalence.
The automorphisms of $\sum\nolimits_{\overline{M}}$ play a key role in this
classification. If $\alpha :\sum\nolimits_{\overline{M}}\rightarrow
\sum\nolimits_{\overline{M}}$ is an automorphism, it is shown in \cite{C2}
that $\alpha $ lifts to an automorphism $\mathbf{R}^{n}\rightarrow \mathbf{R}
^{n},$ which may then be represented by a matrix $A\left( \alpha \right) ,$
as indicated in the following commutative diagram 
\begin{equation*}
\begin{array}{ccc}
\mathbf{R}^{n} & \overset{A\left( \alpha \right) }{\rightarrow } & \mathbf{R}
^{n} \\ 
^{\pi _{\overline{M}}}\downarrow &  & \downarrow ^{\pi _{\overline{M}}} \\ 
\sum_{\overline{M}} & \overset{\alpha }{\rightarrow } & \sum_{\overline{M}}
\end{array}
.
\end{equation*}
Any automorphism $A:\mathbf{R}^{n}\rightarrow \mathbf{R}^{n}$ induces the
map $\widehat{A}:S^{n-1}\rightarrow S^{n-1}$
\begin{equation*}
\widehat{A}\left( x\right) =\dfrac{A\left( x\right) }{\left\| A\left(
x\right) \right\| }.
\end{equation*}
Translating the classification of linear flows into the language of
perspective leads to the following notions.

\begin{definition}
For $X\subset S^{n-1}$ 
\begin{equation*}
\overline{M}\left( X\right) \overset{\text{def}}{=}\cup \left\{ \widehat{
A\left( \alpha \right) }\left( X\right) |\alpha \in Aut\left( \sum\nolimits_{
\overline{M}}\right) \right\} ,
\end{equation*}
and the \emph{rotation class }of a flow $\phi $ on $\sum\nolimits_{\overline{
M}}$ is given by 
\begin{equation*}
\rho \left( \phi \right) \overset{\text{def}}{=}\overline{M}\left( \left\{ 
\text{points of perspective of }\phi \right\} \right) .
\end{equation*}
\end{definition}

\noindent \textbf{Note }Since $x\mapsto -x$ is an automorphism of $
\sum\nolimits_{\overline{M}}$ for any $\overline{M}$, we could work in real
projective $n-1$ space just as well as in $S^{n-1}.$ When applied to the
torus, $\rho \left( \phi \right) $ is called the rotation orbit of $\phi $
in \cite[Chapter 6]{ABZ}.

The topological classification of linear flows given in \cite{C2} may then
be rephrased as follows: 
\begin{equation*}
\left[ \Phi _{\overline{M}}^{\mathbf{\omega }}\text{ is topologically
equivalent to }\Phi _{\overline{M}}^{\mathbf{\omega }^{\prime }}\right] \iff 
\left[ \rho \left( \Phi _{\overline{M}}^{\mathbf{\omega }}\right) =\rho
\left( \Phi _{\overline{M}}^{\mathbf{\omega }^{\prime }}\right) \right] .
\end{equation*}
We now investigate how this can be generalized to other flows.

\begin{lemma}
\label{trans}If $\phi $ and $\psi $ are flows on $\sum\nolimits_{\overline{M}
}$ and $h:\phi \overset{\text{top}}{\approx }\psi $ with $h$ homotopic to
translation by some $\tau \in \sum\nolimits_{\overline{M}}$, then 
\begin{equation*}
\left\{ \text{points of perspective of }\phi \right\} =\left\{ \text{points
of perspective of }\psi \right\} .
\end{equation*}
\end{lemma}

\noindent \textbf{Proof} Suppose $\omega \in S^{n-1}$ is a point of
perspective of $\phi $ and that $\mathcal{O}$ is a lifted $\phi -$orbit of $
x_{0}\in \sum\nolimits_{\overline{M}}$ so that for some sequence of real
numbers $\left\{ t_{i}\right\} _{i\in \mathbf{N}}\rightarrow \infty $ $
\left\{ \mathcal{O}\left( t_{i}\right) \right\} _{i\in \mathbf{N}
}\rightarrow \omega $ in the compactification of perspective. Generally as $
t\rightarrow \infty $ the orbits of the lifted linear flow $\widetilde{\Phi
_{\overline{M}}^{\mathbf{\omega }}}$ approach $\omega $ asymptotically in
the compactification of perspective. Let $\mathbf{R}^{n}\times \left\{
c\right\} $ be the path component of $E_{\overline{M}}$ containing $\mathcal{
O}\left( \mathbf{R}\right) ,\ $which is foliated by $\widetilde{\Phi _{
\overline{M}}^{\mathbf{\omega }}}.$ With each $x\in \mathbf{R}^{n}\times
\left\{ c\right\} $ we associate the ray 
\begin{equation*}
A_{x}\overset{\text{def}}{=}\widetilde{\Phi _{\overline{M}}^{\mathbf{\omega }
}}\left( [0,\infty ),x\right) .
\end{equation*}
Notice that $\overline{A_{x}},$ the closure of $A_{x}$ in the
compactification of perspective, has only $\omega $ in its remainder. Then
for each $i\in \mathbf{N}$ we let $A_{i}$ denote $A_{\mathcal{O}\left(
t_{i}\right) }$, so that $\left\{ \overline{A_{i}}\right\} _{i\in \mathbf{N}
}\rightarrow \left\{ \omega \right\} $ in the Hausdorff metric.

If $\tau \left( \phi \right) $ denotes the flow 
\begin{equation*}
\mathbf{R}\times \sum\nolimits_{\overline{M}}\rightarrow \sum\nolimits_{
\overline{M}};\text{ }\tau \left( \phi \right) \left( t,x\right) =\phi
\left( t,x-\tau \right) +\tau 
\end{equation*}
and if $\mathcal{O}^{\prime }$ is a lifted $\tau \left( \phi \right) -$orbit
of $\tau +x_{0}\in \sum\nolimits_{\overline{M}}$ in $\mathbf{R}^{n}\times
\left\{ c^{\prime }\right\} ,$ then we have $\left( s_{0},c\right) =\mathcal{
O}\left( 0\right) \in \mathbf{R}^{n}\times \left\{ c\right\} $ and $\left(
s_{0}^{\prime },c^{\prime }\right) =\mathcal{O}^{\prime }\left( 0\right) \in 
\mathbf{R}^{n}\times \left\{ c^{\prime }\right\} $ satisfying 
\begin{eqnarray*}
\pi _{\overline{M}}\left( s_{0}\right) +c &=&x_{0}\text{ and }\pi _{
\overline{M}}\left( s_{0}^{\prime }\right) +c^{\prime }=x_{0}+\tau ,\text{
implying} \\
\pi _{\overline{M}}\left( s_{0}^{\prime }-s_{0}\right)  &=&c-c^{\prime
}+\tau ,
\end{eqnarray*}
which in turn implies that if $\widetilde{\tau }\left( s,c\right) \overset{
\text{def}}{=}\left( s+s_{0}^{\prime }-s_{0},c^{\prime }\right) ,$ then 
\begin{equation*}
\begin{array}{ccc}
\mathbf{R}^{n}\times \left\{ c\right\}  & \overset{\widetilde{\tau }}{
\longrightarrow } & \mathbf{R}^{n}\times \left\{ c^{\prime }\right\}  \\ 
\left( \pi _{\overline{M}}+c\right) \downarrow  &  & \downarrow \left( \pi _{
\overline{M}}+c^{\prime }\right)  \\ 
\sum\nolimits_{\overline{M}} & \overset{+\tau }{\longrightarrow } & 
\sum\nolimits_{\overline{M}}
\end{array}
\end{equation*}
is a commutative diagram of maps and $\widetilde{\tau }\ $carries orbits of $
\widetilde{\Phi _{\overline{M}}^{\mathbf{\omega }}}$ to orbits of $
\widetilde{\Phi _{\overline{M}}^{\mathbf{\omega }}}$ and takes $A_{i}$ in $
\mathbf{R}^{n}\times \left\{ c\right\} $ to a corresponding ray $
A_{i}^{\prime }$ of a $\widetilde{\Phi _{\overline{M}}^{\mathbf{\omega }}}-$
orbit in $\mathbf{R}^{n}\times \left\{ c^{\prime }\right\} .$ The rays $
A_{i}^{\prime }$ are translates of the rays $A_{i}$ by $s_{0}^{\prime }-s_{0}
$ and this difference is negligible under the mapping $h_{n}$ as $
i\rightarrow \infty .$ Again taking closures in the compactification of
perspective, $\left\{ \overline{A_{i}^{\prime }}\right\} _{i\in \mathbf{N}
}\rightarrow \left\{ \omega \right\} $ in the Hausdorff metric, implying
that $\left\{ \mathcal{O}^{\prime }\left( t_{i}\right) \right\} _{i\in 
\mathbf{N}}\rightarrow \omega $.

Now the map $\delta \left( x\right) \overset{\text{def}}{=}h\left( x\right)
-\left( \tau +x\right) $ is homotopic to the constant map $\sum\nolimits_{
\overline{M}}\rightarrow \left\{ e_{\overline{M}}\right\} ,$ and so there is
lift $\tilde{\delta}$ making the following diagram commute 
\begin{equation*}
\begin{array}{ccc}
&  & \mathbf{R}^{n} \\ 
& \widetilde{\delta }\nearrow & \downarrow \pi _{\overline{M}} \\ 
\sum\nolimits_{\overline{M}} & \overset{{\large \delta }}{\longrightarrow }
& \sum\nolimits_{\overline{M}}
\end{array}
.
\end{equation*}
Since $\sum\nolimits_{\overline{M}}$ is compact, $\left\| \widetilde{\delta }
\left( x\right) \right\| $ has a maximum value $m.$ Also, with 
\begin{equation*}
\widetilde{h}\left( s,c\right) \overset{\text{def}}{=}\left( s+s_{0}^{\prime
}-s_{0}+\widetilde{\delta }\left( \pi _{\overline{M}}\left( s\right)
+c\right) ,c^{\prime }\right) =``\left( \widetilde{\delta }\left( \pi _{
\overline{M}}\left( s\right) +c\right) ,c^{\prime }\right) +\widetilde{\tau }
\left( s,c\right) ,"
\end{equation*}
the following diagram commutes 
\begin{equation*}
\begin{array}{ccc}
\mathbf{R}^{n}\times \left\{ c\right\} & \overset{\widetilde{h}}{
\longrightarrow } & \mathbf{R}^{n}\times \left\{ c^{\prime }\right\} \\ 
\left( \pi _{\overline{M}}+c\right) \downarrow &  & \downarrow \left( \pi _{
\overline{M}}+c^{\prime }\right) \\ 
\sum\nolimits_{\overline{M}} & \overset{h}{\longrightarrow } & 
\sum\nolimits_{\overline{M}}
\end{array}
.
\end{equation*}
Thus, identifying $\mathbf{R}^{n}\times \left\{ c^{\prime }\right\} $ with $
\mathbf{R}^{n},$ $\left\| \widetilde{h}\left( s,c\right) -\widetilde{\tau }
\left( s,c\right) \right\| \leq m.$ For any preassigned number $K$ we may
choose\ $N$ sufficiently large so that for all $i\geq N$ $\left\{ \inf
\left\| \widetilde{\tau }\left( a\right) \right\| |a\in A_{i}\right\} >K$ in 
$\mathbf{R}^{n}\times \left\{ c^{\prime }\right\} .$ By the construction of $
h_{n}$ and the uniform bound $\left\| \widetilde{h}\left( s,c\right) -
\widetilde{\tau }\left( s,c\right) \right\| \leq m,$ we conclude that $
\left\{ \overline{\widetilde{h}\left( A_{i}\right) }\right\} _{i\in \mathbf{N
}}\rightarrow \left\{ \omega \right\} $ in the compactification of
perspective. The analysis is very similar to that of the torus case (
\cite[Chapter 6,1.5]{ABZ}), where the difference function is not only
uniformly bounded but also doubly periodic. Each $\widetilde{h}\left(
A_{i}\right) $ contains a point $\widetilde{h}\left( \mathcal{O}\left(
t_{i}\right) \right) =\mathcal{O}^{\prime \prime }\left( t_{i}^{\prime
\prime }\right) $ in a lift $\mathcal{O}^{\prime \prime }$ of the $\psi -$
orbit of $h\left( x_{0}\right) .$ While there is no reason to expect that $
t_{i}^{\prime \prime }=t_{i},$ we may conclude $\left\{ t_{i}^{\prime \prime
}\right\} _{i\in \mathbf{N}}\rightarrow \infty $ since $h$ preserves the
orientation of orbits and $\left\{ t_{i}\right\} _{i\in \mathbf{N}
}\rightarrow \infty .$ Hence, $\omega $ is a point of perspective of $\psi $
as well. The other set containment follows mutatis mutandis. \hfill $\square 
$

\noindent \textbf{Note} This implies that $\rho \left( \phi \right) =\rho
\left( \psi \right) .$ It follows by setting $\tau =e_{\overline{M}}$ in the
above proof that if $\mathcal{O}$ is a lifted orbit having $\omega $ as a
forward limit point, then any other lift of the same orbit also has $\omega $
as a forward limit point.

\begin{lemma}
If $\alpha \in Aut\left( \sum\nolimits_{\overline{M}}\right) $ and $\alpha
:\phi \overset{\text{top}}{\approx }\psi ,$ then $\rho \left( \phi \right)
=\rho \left( \psi \right) .$
\end{lemma}

\noindent \textbf{Proof} Let $\omega \in S^{n-1}$, $\mathcal{O}$ , $x_{0}\in
\sum\nolimits_{\overline{M}}$ and $\left\{ t_{i}\right\} _{i\in \mathbf{N}
}\rightarrow \infty $ be as above. With $\omega ^{\prime }=\widehat{A\left(
\alpha \right) }\left( \omega \right) ,$ $\alpha $ maps orbits of $\Phi _{
\overline{M}}^{\mathbf{\omega }}$ to orbits of $\Phi _{\overline{M}}^{
\mathbf{\omega }^{\prime }}$ after (possibly) rescaling time (see \cite[3.5]
{C2}), and $\alpha $ lifts to an affine map $\mathbf{R}^{n}\times \left\{
c\right\} \rightarrow \mathbf{R}^{n}\times \left\{ c^{\prime }\right\} $, so
that the above argument for $\tau \left( \phi \right) $ can be slightly
adapted and applied here to show that $\omega ^{\prime }$ is a point of
perspective of $\psi .$ Since $\overline{M}\left( \left\{ \omega \right\}
\right) =\overline{M}\left( \left\{ \omega ^{\prime }\right\} \right) ,$ we
conclude that $\rho \left( \phi \right) \subset \rho \left( \psi \right) .$
Similarly, with $\alpha ^{-1}$ in place of $\alpha ,$ $\rho \left( \phi
\right) \supset \rho \left( \psi \right) .$ \hfill $\square $

\begin{theorem}
If the flows $\phi $ and $\psi $ on $\sum\nolimits_{\overline{M}}$ are
topologically equivalent, then $\rho \left( \phi \right) =\rho \left( \psi
\right) .$
\end{theorem}

\noindent \textbf{Proof }Suppose that $h:\phi \overset{\text{top}}{\approx }
\psi $ and let $\omega \in S^{n-1}$, $\mathcal{O}$ , $x_{0}\in
\sum\nolimits_{\overline{M}}$ and $\left\{ t_{i}\right\} _{i\in \mathbf{N}
}\rightarrow \infty $ be as above. With $\tau :\sum\nolimits_{\overline{M}
}\rightarrow \sum\nolimits_{\overline{M}}$ denoting translation by $h\left(
e_{\overline{M}}\right) ,$ the homeomorphism $g=\tau ^{-1}\circ h$ fixes the
identity element $e_{\overline{M}}.$ From \cite{S} (or see \cite[3.8]{C2}),
it follows that there is a homotopy
\begin{equation*}
H:\left[ 0,1\right] \times \sum\nolimits_{\overline{M}}\rightarrow
\sum\nolimits_{\overline{M}}
\end{equation*}
from $g=H_{0}$ to an automorphism $\alpha =H_{1}.$ Then 
\begin{equation*}
h^{\prime }\overset{\text{def}}{=}h\circ \alpha ^{-1}=\tau \circ g\circ
\alpha ^{-1}
\end{equation*}
is homotopic to $\tau $ via $\tau \circ H_{t}\circ \alpha ^{-1},$ implying
that $h=h^{\prime }\circ \alpha ,$ where $h^{\prime }$ is homotopic to a
translation. Then by the above lemmas it follows that $\omega ^{\prime }=
\widehat{A\left( \alpha \right) }\left( \omega \right) $ is a point of
perspective of $\psi .$ Since $\overline{M}\left( \left\{ \omega \right\}
\right) =\overline{M}\left( \left\{ \omega ^{\prime }\right\} \right) ,$ we
conclude that $\rho \left( \phi \right) \subset \rho \left( \psi \right) $
and similarly $\rho \left( \phi \right) \supset \rho \left( \psi \right) .$
\hfill $\square $

Notice that if we replace $\left\{ t_{i}\right\} _{i\in \mathbf{N}
}\rightarrow \infty $ with $\left\{ t_{i}\right\} _{i\in \mathbf{N}
}\rightarrow \pm \infty $ in definition \ref{persp}, we could obtain a
similar invariant. Also, in analogy with the invariant introduced in \cite
{MZ}, one could consider all limit points of sequences of the form $\left\{ 
\mathcal{O}_{x_{i}}\left( t_{i}\right) \right\} _{i\in \mathbf{N}}$ where $
x_{i}$ is allowed to vary and  $\left\{ t_{i}\right\} _{i\in \mathbf{N}
}\rightarrow \infty .$

\section{Exponents and Other Invariants}

We now compare $\rho \left( \phi \right) $ with the exponent group of an
orbit $\phi _{x}:\mathbf{R}\rightarrow X$ of a flow $\phi $ as introduced in 
\cite{C1}. The exponent group of $\phi _{x}$\ is the subgroup of $\left( 
\mathbf{R},+\right) $ given by 
\begin{equation*}
\left\{ \alpha \in \mathbf{R\mid }
\begin{array}{c}
\left\{ \pi ^{1}\left( \alpha t_{i}\right) \right\} \text{ converges in }
S^{1}\text{ for all sequences }\left\{ t_{i}\right\} \text{ } \\ 
\text{for which }\left\{ \phi _{x}\left( t_{i}\right) \right\} \text{
converges in }X
\end{array}
\right\} .
\end{equation*}
Recall that a non-empty set $\mathcal{M}\subset X$ is a \emph{minimal set}
of the flow $\phi $ on $X$ if and only if $\mathcal{M}$ is closed and
invariant and has no proper non-empty subset which is also closed and
invariant. In \cite[Theorem 5]{C1} the exponent group is shown to be the
same for all the orbits of a minimal set of a flow, allowing one to
associate a group $\mathcal{E}\left( \mathcal{M}\right) $ with each minimal
set $\mathcal{M}$. For example, an $\alpha -$irrational flow on the torus
formed by taking the suspension of rotation by $\alpha $ on $S^{1}$ is
minimal and its exponent group is $\left\{ m+n\alpha \mathbf{\mid }m,n\in 
\mathbf{Z}\right\} .$ While this group is a $C^{0}$-conjugacy invariant 
\cite[Corollary 6]{C1}, it is not invariant under time changes since there
are irrational flows which are not conjugate to flows to which they are
topologically equivalent (see, e.g., \cite[V.8.19]{NS}).\ Such time changed
flows will have the same rotation classes as the original flows but
different exponent groups. The exponent group is sensitive to changes in the
ergodic properties of the flows (stability) which the rotation class does
not detect. Which is more desirable depends on what one wishes to
investigate: topological equivalence or ergodic properties.

If $h:X\rightarrow Y$ provides a continuous semiconjugacy between a flow $
\phi $ on $X$ and a flow $\psi $ on $Y$, $\phi $ is said to be an \emph{
almost 1:1-extension }of $\psi $ when there is some $x\in X$ with $
h^{-1}\left( h\left( x\right) \right) =\left\{ x\right\} ,$ see, e.g., 
\cite[Chapter IV,6.1]{V}. The following theorem allows one to calculate $
\mathcal{E}\left( \mathcal{M}\right) $ for a large class of minimal sets.
Recall that an orbit of a flow $\phi $ on a compact metric space $X$ is 
\emph{almost periodic }if and only if it is uniformly Lyapunov stable (i.e.,
equicontinuous: see, e.g., \cite[Chapter V.11.7]{NS}).

\begin{theorem}
If the minimal flow $\phi $ on a compact space $\mathcal{M}$ is an almost
1:1-extension of the almost periodic flow $\psi $ on $\Sigma ,$ then $
\mathcal{E}\left( \mathcal{M}\right) =\mathcal{E}\left( \Sigma \right) .$
\end{theorem}

\noindent \textbf{Proof }This argument depends critically on the results and
constructions to be found in \cite{C1}. We assume that the flow $\psi $ is
represented as a linear flow on a $\kappa -$solenoid $\Sigma $ with identity 
$e$. As $\psi $ is invariant under translations, we may assume that the
semiconjugacy $h:\mathcal{M}\rightarrow \Sigma $ is such that $h^{-1}\left(
h\left( x\right) \right) =h^{-1}\left( e\right) =\left\{ x\right\} $ for
some $x\in \mathcal{M}$. From \cite[Theorem 5]{C1} it follows that $\mathcal{
E}\left( \mathcal{M}\right) \supset \mathcal{E}\left( \Sigma \right) $.
Generally, $\mathcal{E}\left( \mathcal{M}\right) $\ coincides with the
exponent group of a maximal irrational linear flow $\mu $ on a topological
group $\mathcal{G}$ semiconjugate to $\phi ,$ and $\mathcal{E}\left( 
\mathcal{M}\right) $ is naturally isomorphic with the dual group of the
topological group supporting the flow $\mu $, which is isomorphic to $\check{
H}^{1}\left( \mathcal{G}\right) $. If $g:\mathcal{M}\rightarrow \mathcal{G}$
provides a semiconjugacy of $\phi $ with $\mu $, again we may assume that $
g\left( x\right) $ is the identity of $\mathcal{G}$. If the containment $
\mathcal{E}\left( \Sigma \right) \subset \mathcal{E}\left( \mathcal{M}
\right) $ is proper, the Pontryagin dual epimorphism $E:\mathcal{
G\rightarrow }\Sigma $ is not an isomorphism. Then $E$ also semiconjugates
the flow $\mu $ with the flow $\psi $. Since the flows involved are minimal,
two semiconjugacies coinciding at a point must be identical. From this it
follows that $E\circ g=h$. But $E$ is not 1:1 and $g$ is surjective,
contradicting the assumption that $h^{-1}\left( h\left( x\right) \right) $ $
=\left\{ x\right\} .$ We must therefore have $\mathcal{E}\left( \Sigma
\right) =\mathcal{E}\left( \mathcal{M}\right) .$ \hfill $\square $

Since the natural semiconjugacy of a Denjoy flow with the irrational flow
from which it is constructed is an almost 1:1 extension, this allows us to
calculate its exponent group easily. For the natural flow on an $\alpha -$
Denjoy minimal set $\mathbb{D}_{\alpha }$ (regardless of the number of
components in the complement) the exponent group is $\left\{ m+n\alpha 
\mathbf{\mid }m,n\in \mathbf{Z}\right\} .$ Translating Fokkink's
classification of the Denjoy continua with connected complements \cite{F},
\cite{BW}, we have: 
\begin{gather*}
\left[ \mathbb{D}_{\alpha }\text{ is homeomorphic to }\mathbb{D}_{\beta }
\right]  \\
\Updownarrow  \\
\left[ \mathcal{E}\left( \mathbb{D}_{\alpha }\right) =c\cdot \mathcal{E}
\left( \mathbb{D}_{\beta }\right) \text{ for some constant }c\right]  \\
\Updownarrow  \\
\left[ \overline{\left( 
\begin{array}{cc}
1 & 0 \\ 
0 & 1
\end{array}
\right) }\left( \dfrac{\left( 1,\alpha \right) }{\left\| \left( 1,\alpha
\right) \right\| }\right) =\overline{\left( 
\begin{array}{cc}
1 & 0 \\ 
0 & 1
\end{array}
\right) }\left( \dfrac{\left( 1,\beta \right) }{\left\| \left( 1,\beta
\right) \right\| }\right) \right] 
\end{gather*}
which is equivalent to requiring that $\alpha $ and $\beta $ have continued
fraction expansions with a common tail \cite{BW},\cite[Chapter 6.1.7]{ABZ},
and so all of the invariants under consideration suffice to classify these
minimal sets \emph{topologically}. (Recall that 
\begin{equation*}
\left\{ \widehat{A\left( \alpha \right) }|\alpha \in Aut\left( \mathbf{T}
^{2}\right) \right\} =GL\left( 2,\mathbf{Z}\right) .)
\end{equation*}

It is then natural to wonder whether there is a 1-dimensional compact
minimal flow having a given non-trivial countable subgroup $\mathcal{E}
\subset \left( \mathbf{R,+}\right) $ as its exponent group. We now construct
such a flow for subgroups with finite torsion-free rank (see  \cite{Fu}).
First we find a maximal rationally independent subset 
\begin{equation*}
A=\left\{ \alpha _{i}\right\} _{i=1}^{n}\subset \mathcal{E}
\end{equation*}
($n$ is the torsion-free rank of $\mathcal{E}$) with 
\begin{equation*}
\left\{ \alpha _{i}\right\} _{i=1}^{n-1}\subset \mathbf{R-Q}\text{ }(\text{
it is possible that }A\subset \mathbf{R-Q}).
\end{equation*}
If $n=1$ we obtain the desired flow by a linear flow on a $1-$solenoid. When 
$n>1$, for each $i\in \left\{ 1,..,n-1\right\} $ we construct a Denjoy
homeomorphism $h_{i}$ of $S^{1}$ with rotation number $\alpha _{i}$. Then
for each $i$ there is a monotone map $\sigma _{i}:S^{1}\rightarrow S^{1}$
isotopic to $id_{S^{1}}$ which provides a semiconjugacy between $h_{i}$ and
rotation by $\pi ^{1}\left( \alpha _{i}\right) $ in $S^{1}.$ Each such $h_{i}
$ has a minimal Cantor Set $C_{i}$. Then 
\begin{equation*}
C=\prod_{i=1}^{n-1}C_{i}
\end{equation*}
is yet another Cantor Set and 
\begin{equation*}
h\left( \left( x_{i}\right) _{i=1}^{n-1}\right) \overset{\text{def}}{=}
\left( h_{i}\left( x_{i}\right) \right) _{i=1}^{n-1}
\end{equation*}
is a homeomorphism of $\mathbf{T}^{n-1}$ that has $C$ as an invariant set
since each $C_{i}$ is $h_{i}-$invariant. Also 
\begin{equation*}
\sigma \left( \left( x_{i}\right) _{i=1}^{n-1}\right) \overset{\text{def}}{=}
\left( \sigma _{i}\left( x_{i}\right) \right) _{i=1}^{n-1}
\end{equation*}
is isotopic to $id_{\mathbf{T}^{n-1}}$ and provides a semiconjugacy between $
h$ and $R,$ rotation in $\mathbf{T}^{n-1}$ by $\left( \pi ^{1}\left( \alpha
_{i}\right) \right) _{i=1}^{n-1},$ which is minimal by the rational
independence of the $\alpha _{i}.$ The semiconjugacy of $h|_{C}$ to $R$
shows that $h|_{C}$ is minimal \cite[Chapter IV.6.1(b)]{V} and $h|_{C}$ is
an almost 1:1 extension of $R.$ The suspension of $h$ will then be a flow on 
$\mathbf{T}^{n}$ having a one--dimensional minimal set having $C$ as a $0-$
dimensional cross-section. By adjusting the time scale on this flow to have
a return time to $C$ of $1/\alpha _{n}$, we obtain a minimal flow $\phi $ on
a one--dimensional subcontinuum of $\mathbf{T}^{n},$ and the orbits of this
flow have as their exponent group $\left\langle A\right\rangle $, the
subgroup of $\left( \mathbf{R},+\right) $ generated by $A$ since this flow
is an almost 1:1 extension of the linear flow $\Phi ^{\mathbf{\alpha }}$ on $
\mathbf{T}^{n}$ having exponent group $\left\langle A\right\rangle $ one
obtains by taking the suspension of $R$ and rescaling time by a factor of $
1/\alpha _{n}.$

By enumerating the elements of $\mathcal{E}-\left\langle A\right\rangle
=\left\{ g_{1},g_{2},....\right\} $ (possibly empty), we obtain a direct
limit representation of $\mathcal{E}$
\begin{equation*}
\left\langle A\right\rangle \hookrightarrow \left\langle A\cup \left\{
g_{1}\right\} \right\rangle \hookrightarrow \left\langle A\cup \left\{
g_{1},g_{2}\right\} \right\rangle \hookrightarrow \cdots \mathcal{E}
\end{equation*}
and dual to this direct sequence is an inverse sequence 
\begin{equation*}
\mathbf{T}^{n}\overset{E_{1}}{\twoheadleftarrow }\mathbf{T}^{n}\overset{E_{2}
}{\twoheadleftarrow }\mathbf{T}^{n}\overset{E_{3}}{\twoheadleftarrow }\cdots
\sum\nolimits_{\overline{M}}
\end{equation*}
where each $E_{i}$ is an epimorphism and the inverse limit an $n-$solenoid.
Pulling the flow $\phi $ back by the various $E_{i}$ as described in \cite
{C3}, we obtain a flow on $\sum\nolimits_{\overline{M}}$. This flow has a
one--dimensional invariant set $\mathcal{M}$ which is minimal \cite[Chapter
IV.6.1(b)]{V} and is an almost 1:1 extension of the corresponding linear
flow obtained by pulling back $\Phi ^{\mathbf{\alpha }}$ to a linear flow on 
$\sum\nolimits_{\overline{M}}.$ The orbits of the flow on $\mathcal{M}$ have 
$\mathcal{E}$ as their exponent group since this is the exponent group of
the corresponding linear flow on $\sum\nolimits_{\overline{M}},$ see 
\cite[Section 4]{C1}. This construction in the case $n=2$ is considered in
detail in \cite{C3} where such minimal sets (called denjoids) are classified
topologically. In the case of a general denjoid the exponent group is not
finitely generated and has torsion-free rank 2. 

These flows will have only one point of perspective:
\begin{equation*}
\mathbf{\alpha =}\dfrac{\left( \alpha _{1},...,\alpha _{n}\right) }{\left\|
\left( \alpha _{1},...,\alpha _{n}\right) \right\| }
\end{equation*}
since the lifted flow lines in $\mathbf{R}^{n}\times \left\{ c\right\} $
differ from a lifted irrational flow by a map (obtained from $\sigma \left(
x\right) -x)$ homotopic to a constant, implying the existence of a bounded $
\tilde{\delta}$ as in \ref{trans} (see \cite[Chapter 6.1.7]{ABZ} for similar
arguments). Once again we find that the exponent group and the rotation
class give us the same information when analyzing these flows on a fixed $
\sum\nolimits_{\overline{M}}$ since the rotation classes of such flows are
identical precisely when their exponent groups are multiples of each other.
However, the exponent group is independent of the embedding in $
\sum\nolimits_{\overline{M}}$ and so can be used to compare flows on $
\sum\nolimits_{\overline{M}}$ and $\sum\nolimits_{\overline{N}}$ quite
easily. To translate information about the rotation class of a flow on $
\sum\nolimits_{\overline{M}}$ to information about the class of a flow on $
\sum\nolimits_{\overline{N}}$ when $\sum\nolimits_{\overline{M}}$ and $
\sum\nolimits_{\overline{N}}$ are isomorphic, one can use a matrix
representing an isomorphism $i:\sum\nolimits_{\overline{M}}\rightarrow
\sum\nolimits_{\overline{N}},$ see \cite[Theorem 3.4]{C2}.

We now compare the rotation class with the rotation sets as discussed in 
\cite{FM}. There examples are discussed with rotation sets that are
intervals. In particular, the example due to Katok of a time changed
irrational flow with a singular point has an interval as its rotation set,
but all the flow lines with unbounded lifted forward orbits would have the
same point of perspective as the original irrational flow. Once again, which
invariant is best depends on the goal of the investigation.

\end{document}